\documentclass[
nopreprint,%
showpacs,
showkeys
]{jasatex}

\usepackage{graphicx}
\usepackage{dcolumn}
\usepackage{amsmath,amsfonts}
\usepackage{bm}
\usepackage{subfigure}
\usepackage{color}
\usepackage[normalem]{ulem}
\usepackage{pgf,tikz}
\usetikzlibrary{arrows}

\usepackage[show]{ed} %

\newcommand{\rf}[1]{\eqref{#1}}
\newcommand{\He}{{\rm H}}

\newcommand{\ri}{{\rm i}}
\newcommand{\re}{{\rm e}}
\newcommand{\rd}{{\rm d}}
\DeclareMathOperator{\sech}{sech}
\DeclareMathOperator{\gd}{gd}

\begin{document}

\title[Diffraction by an impedance wedge]{Diffraction by a right-angled impedance wedge: an edge source formulation}

\author{David P. Hewett\footnote{Author to whom correspondence should be addressed. \\ Electronic mail: \textsf{hewett@maths.ox.ac.uk}}}
\affiliation{Mathematical Institute, University of Oxford, Radcliffe Observatory
Quarter, Woodstock Road, Oxford, OX2 6GG, United Kingdom}
\author{Aaron Morris}
\affiliation{Mathematical Institute, University of Oxford, Radcliffe Observatory
Quarter, Woodstock Road, Oxford, OX2 6GG, United Kingdom}

\date{\today}

\begin{abstract}

This paper concerns the frequency domain problem of diffraction of a plane wave incident on an infinite right-angled wedge on which impedance (absorbing) boundary conditions are imposed. It is demonstrated that the exact Sommerfeld-Malyuzhinets contour integral solution for the diffracted field can be transformed to a line integral over a physical variable along the diffracting edge. This integral can be interpreted as a superposition of secondary point sources (with directivity) positioned along the edge, in the spirit of the edge source formulations for rigid (sound-hard) wedges derived in [U.~P.~Svensson, P.~T.~Calamia and S.~Nakanishi, Acta Acustica/Acustica 95, 2009, pp.~568-572]. However, when surface waves are present the physical interpretation of the edge source integral must be altered: it no longer represents solely the diffracted field, but rather includes surface wave contributions.
\end{abstract}


\maketitle

\section{Introduction}\label{sect:intro}

Diffraction by an infinite wedge is a fundamental canonical problem in acoustic scattering. Exact closed-form frequency-domain solutions for point source, line source or plane wave excitation with homogeneous Dirichlet (sound soft) or Neumann (sound hard, or rigid) boundary conditions are available in many different forms\cite{Macdonald,BowmanSenior,Pierce}. For example, series expansions in terms of eigenfunctions are available for near field calculations (e.g.~for analysing edge singularities). Contour integral representations over so-called Sommerfeld-Malyuzhinets contours are better suited to far field computations (e.g.~for deriving diffraction coefficients in computational methods such as the Geometrical Theory of Diffraction \cite{Keller62,Borovikov}). 
More recently it has been discovered that the `diffracted' component of these solutions (precisely, that which remains after subtracting from the total field the geometrical acoustics terms) can be expressed in a more physically intuitive form, namely as a line integral superposition of directional secondary sources located along the diffracting edge \cite{SvenssonCalamia2009}.  (In fact the frequency domain expressions derived in Ref.~\cite{SvenssonCalamia2009} had appeared already in Ref.~\cite{Buckingham1989}, but the interpretation in terms of secondary edge sources seems to have been first made in Ref.~\cite{SvenssonCalamia2009}.) 

One appealing feature of the edge source interpretation is that it offers a natural way to write down approximate solutions for finite edges, simply by truncating the domain of integration. It has also led to edge integral equation formulations of scattering problems, where the integral equation is posed on the union of all the scatterer's edges \cite{AsheimSvensson13}. For Dirichlet and Neumann boundary conditions the edge source formulations are now well understood: efficient numerical evaluation of the line integrals has been considered in Ref.~\cite{AsheimSvensson10} using the method of numerical steepest descent, as has the behaviour of the line integrals near shadow boundaries \cite{SvenssonCalamia2006} and edges \cite{HewettSvensson13}. Note also Ref.~\cite{Svensson1999}, where the corresponding time-domain case is considered.

For the more difficult case of diffraction by a wedge with impedance (absorbing) boundary conditions, some exact solutions are also known. For example, the case of plane wave incidence on an impedance wedge can be solved using the Sommerfeld-Malyuzhinets technique, and converted to a series expansion using a Watson-type transformation (see Ref.~\cite{Osipov99} and the references therein). But the solution obtained is much more cumbersome than those for the corresponding Dirichlet and Neumann problems, and the technique requires the solution of a certain non-trivial functional difference equation. This increased complexity is perhaps to be expected, since the physics of the impedance problem are fundamentally more complicated than those for the Dirichlet and Neumann problems; in particular, the wedge faces can under certain conditions support surface waves.

However, for the special (yet important) case of a right-angled wedge, the solution takes a particularly simple and explicit form\cite{Rawlins90}. In Ref.~\cite{Rawlins90}, Rawlins proves that the solution to the impedance problem for a right-angled wedge (with possibly different impedances on each face) can be obtained from that of the corresponding Dirichlet problem, generalised to allow complex incident angles, by the application of a certain linear differential operator (see Eqs.~\rf{eqn:pIRep}--\rf{eqn:tildpDef} below for details). Rawlins applies this operator to the classical series and integral representations of the Dirichlet solution to obtain relatively simple series and integral representations for the impedance solution. (The solution for the case where the impedance is the same on both faces was presented previously in a similar but more complicated form in Ref.~\cite{Pierce78}.)

In this paper it will be shown that Rawlins' solution for the impedance wedge can be transformed into an edge source representation of the same form as those derived for rigid (sound-hard) wedges in Ref.~\cite{SvenssonCalamia2009}. 
This appears to be the first edge source representation for diffraction by an impedance wedge. While the solution obtained is valid only for a right-angled wedge, it should be remarked that this special case is ubiquitous in many acoustical applications (e.g., urban acoustics \cite{Hewett12}). 

\section{Edge source integral for ideal wedges}\label{sect:EdgeSource}

The edge source formulation for ideal (Dirichlet and Neumann) wedges will briefly be reviewed. For the most general setting (illustrated in Fig.~1), consider a point source $\mathrm{S}$ and point receiver $\mathrm{R}$ in the presence of a wedge of exterior angle $\pi<\theta_W<2\pi$. 

Let $(r,\theta,z)$ denote cylindrical coordinates with the $z$-axis along the edge, the propagation domain occupying the region $0<\theta<\theta_W$, and the wedge the region $\theta_W<\theta<2\pi$. Consider also Cartesian coordinates $(x,y,z)$ with $r=\sqrt{x^2+y^2}$, $x=r\cos\theta$, $y=r\sin\theta$. Without loss of generality it will be assumed that the receiver is located in the plane $z=0$, at $\mathrm{R}=(r,\theta,0)$. For each point $z$ on the edge one can also introduce local spherical coordinates $(l,\theta,\phi)$, with $\theta$ defined as before and $l=\sqrt{r^2+z^2}$, $r=l\sin\phi$, $z=-l\cos\phi$.

For consistency with Ref.~\cite{Rawlins90} the time-dependence $\re^{-\ri \omega t}$ will be assumed throughout. 
Then the diffracted field at $\mathrm{R}$ (i.e.\ the total field minus the geometrical acoustics field) due to a monopole source at $\mathrm{S}=(r_0,\theta_0,z_0)$ 
can be written as a line integral 
\begin{equation}
p_{d} = -\frac{\nu}{4 \pi } {{\int}_{\!\!\!\!-\infty}^{\infty}}
\frac {{\rm e}^{\ri k(l_0+l)}} {l_0 l}
\beta\, \rd z
\label{eq:basicintegral}
\end{equation}
over edge positions $z$, where $1/2<\nu=\pi/\theta_W<1$ is the wedge index \cite{SvenssonCalamia2009}. The integral in Eq.~\rf{eq:basicintegral} can be interpreted as a superposition of secondary edge sources along the edge. The factor $\beta$ can be interpreted as a directivity function, and takes the following forms for Dirichlet ($p=0$) and Neumann ($\partial p/\partial n=0$) boundary conditions: %
\[
\beta^{\rm D} = -\beta_1 + \beta_2 + \beta_3 - \beta_4, \quad
\beta^{\rm N} = \beta_1 + \beta_2 + \beta_3 + \beta_4,
\]
where
\[
\beta_i = \frac{\sin \left( \nu \varphi_i   \right)}{\cosh \left( \nu \eta  \right) - \cos \left( \nu \varphi_i  \right)},
\quad i=1,\ldots,4,
\]
\begin{align*}
\varphi_1 = \pi + \theta_0 + \theta, \qquad  \varphi_2 = \pi + \theta_0 - \theta, \\  
\varphi_3 = \pi - \theta_0 + \theta, \qquad  \varphi_4 = \pi - \theta_0 - \theta,   
\end{align*}
and the auxiliary function $\eta$ is
\begin{align}
\label{eq:etaDef}
\eta &= \cosh^{-1}\frac{(z-z_0)z + l_0 l}{r_0 r} %
=\cosh^{-1}\frac{1+\cos \phi_0 \cos \phi}{\sin \phi_0 \sin \phi}.
\end{align}
The second expression for $\eta$ in Eq.~\rf{eq:etaDef} shows that $\beta$ is a function only of the local spherical angles $\theta,\phi$; this justifies the interpretation of $\beta$ as a directivity function. 
The formula for $\beta^D$ above corrects that in Ref.~\cite{HewettSvensson13} by reversing the sign of $\beta^D$. Also, the first expression for $\eta$ in Eq.~(\ref{eq:etaDef}) corrects a sign error in the corresponding formula in Ref.~\cite{SvenssonCalamia2009}.
A mixed wedge with Dirichlet conditions on $\theta=0$ and Neumann conditions on $\theta=\theta_W$ can also be treated, by summing the Dirichlet solutions for a wedge angle $2\theta_W$ and source positions $\theta_0$ and $2\theta_W-\theta_0$ respectively. This gives the directivity factor
\begin{align*}
\label{}
\beta^{DN}= -\beta'_1 + \beta'_2 + \beta'_3 - \beta'_4,
\end{align*}
where 
\begin{align}
\label{eqn:betaiDef}
\beta'_i = \frac{2\sin \left(\frac{\nu}{2} \varphi_i   \right)\cosh \left( \frac{\nu}{2}\eta  \right)}{\cosh{\nu \eta} - \cos \left( \nu\varphi_i  \right)},
\end{align}
which, when inserted in Eq.~\rf{eq:basicintegral}, agrees with the solution derived by Buckingham in Eqs.~(62)-(63) of Ref.~\cite{Buckingham1989} using a modal expansion. 
The solution for a mixed wedge with Neumann conditions on $\theta=0$ and Dirichlet conditions on $\theta=\theta_W$ can then be obtained by replacing $\theta$ by $\theta_W-\theta$ and $\theta_0$ by $\theta_W-\theta_0$ in Eq.~\rf{eqn:betaiDef}.

\begin{figure}[t]
\begin{center}
\includegraphics[width=8cm]{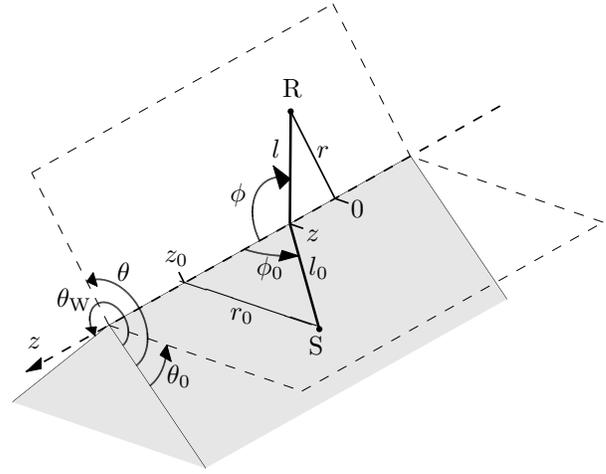}
\caption{The geometry of the wedge.}\label{fig:wedge}
\end{center}
\end{figure}

\subsection{Plane wave incidence perpendicular to the edge}
For the 
two-dimensional case of plane wave incidence perpendicular to the edge, the source is placed at $z_0=0$ with $\phi_0 = \pi/2$. In this case 
\begin{align*}
 \eta = \cosh^{-1}\!  \left(\frac{l}{r} \right)
 \! =  \cosh^{-1}\!  \left(\!\!  \sqrt{1+\frac{z^2}{r^2}} \right)
\! =  \cosh^{-1}\!  \left(\frac{1}{\sin \phi} \right)\! ,
\end{align*}
and, because of symmetry, the integration range in Eq.~(\ref{eq:basicintegral}) can be halved. Furthermore, for plane wave incidence the $1/l_0$ spherical attenuation factor is removed, and, if one refers the phase of the diffracted sound pressure relative to the arrival time at the edge, the phase oscillation factor ${\rm e}^{\ri kl_0}$ also disappears, giving
\begin{align}
p_{d} &= -\frac{\nu}{2 \pi }   \int_{0}^{\infty}
\frac {\ {\rm e}^{\ri kl}} {l}\beta \, \rd z.
\label{eq:pwintegral}
\end{align}

\subsection{Right-angled wedge}%

For plane wave incidence perpendicular to a right-angled wedge, for which $\theta_W=3\pi/2$ (i.e.\ $\nu=2/3$), further simplification is obtained. 
Now Eq.~\rf{eq:pwintegral} is
\begin{align}
p_{d} &= -\frac{1}{3 \pi }   \int_{0}^{\infty}
\frac {\ {\rm e}^{\ri kl}} {l}\beta \, \rd z,
\label{eq:pwintegralspecial}
\end{align}
and using multiple angle formulas one can show that
\begin{align*}
\hspace{-3mm}\cosh \left(\frac{2\eta}{3} \right)
= \frac{1}{2}\left(\lambda^{-2/3}+\lambda^{2/3}\right),
\end{align*}
where
\begin{align*}
\label{}
\lambda = w+\sqrt{w^2-1}, \qquad w=\frac{l}{r} =\sqrt{1+\frac{z^2}{r^2}} = \frac{1}{\sin \phi}.
\end{align*}
To summarise, for plane wave incidence perpendicular to the edge of a right-angled wedge, with $0\leq \theta_0\leq  3\pi/2$, the total field is 
\begin{align}
p (r,\theta)
&=\He[\pi-|\theta-\theta_0|]\re^{-\ri kr\cos(\theta-\theta_0)}\notag \\
&\quad+ R_1\He[\pi-\theta-\theta_0]\re^{-\ri kr\cos(\theta+\theta_0)}\notag \\
&\quad+ R_2\He[\theta+\theta_0-2\pi]\re^{\ri kr\cos(\theta+\theta_0)}\notag \\
&\quad+ p_d(r,\theta),
\label{eqn:ptotal}
\end{align}
where the diffracted field $p_d$ is given by Eq.~\rf{eq:pwintegralspecial} and the remaining three terms represent the geometrical acoustics field (here $\He[t]$ is the Heaviside step function: $\He[t]=1$, $t>0$; $\He[t]=1/2$, $t=0$; $\He[t]=0$, $t<0$). The first term in Eq.~\rf{eqn:ptotal} represents the incident field, the second the reflected wave from the face $\theta=0$ and the third the reflected wave from $\theta=\theta_W$. $R_1,R_2$ are reflection coefficients: for the pure Dirichlet case $R_1=R_2=-1$; for the pure Neumann case $R_1=R_2=1$; and for the mixed Dirichlet/Neumann and Neumann/Dirichlet cases $R_1=-1$, $R_2=1$ and $R_1=1$, $R_2=-1$, respectively.

\subsection{Derivation and generalization to complex incident angles}
\label{sec:DirDerivation}

In order to understand how Eq.~\rf{eqn:ptotal} should be generalised to the case of impedance boundary conditions, it is instructive to review its derivation from the classical Sommerfeld contour integral solution. Part of this derivation was presented already in Ref.~\cite{SvenssonCalamia2009}, but in this section the analysis of Ref.~\cite{SvenssonCalamia2009} will be extended to allow a complex incident angle $\theta_0$, corresponding physically to diffraction by an inhomogeneous plane wave. 
For brevity, attention will be restricted to the Dirichlet case, but the Neumann and mixed problems can be analysed similarly. 

\begin{figure}[t]
\begin{center}
\includegraphics[scale=1]{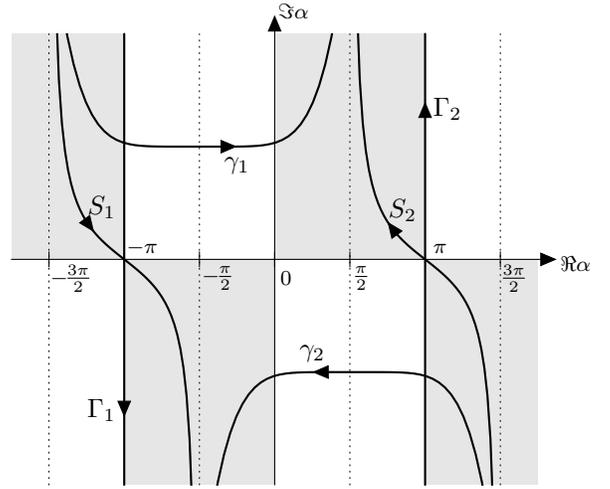}
\caption{Integration contours for Eq.~\rf{eqn:DirSomm}. The shaded sectors are those in which the integrand decays exponentially at infinity.\label{fig:gamma}}
\end{center}
\end{figure}

For the Dirichlet case the Sommerfeld contour integral solution is \cite{BowmanSenior}
\begin{align}
\label{eqn:DirSomm}
p^D(r,\theta) = -\frac{1}{3\pi \ri} \int_{\gamma_1+\gamma_2} \frac{\re^{-\ri k r \cos{\alpha}}\sin \frac{2}{3}\theta_0}{\cos \frac{2}{3}\theta_0 - \cos\frac{2}{3}(\alpha+\theta)}\,\rd \alpha,
\end{align}
where the integration contour is in two parts: $\gamma_1$ lies above all singularities of the integrand and goes from $-\frac{3\pi}{2}+\ri \infty$ to $\frac{\pi}{2}+\ri \infty$, and $\gamma_2$ lies below all singularities of the integrand and goes from $\frac{3\pi}{2}-\ri \infty$ to $-\frac{\pi}{2}-\ri \infty$ (see Fig.~\ref{fig:gamma}). The integral converges rapidly for all complex $\theta$ and $\theta_0$, and the integrand has poles at $\alpha = -\theta\pm \theta_0 + 3n\pi$ for $n\in\mathbb{Z}$. 
As is pointed out in Ref.~\cite{Rawlins90}, the formula in Eq.~\rf{eqn:DirSomm} makes sense not just for real $0\leq \theta_0\leq \frac{3\pi}{2}$, but also for all complex $\theta_0$. 

\begin{figure}[t]
\begin{center}
\subfigure[\label{fig:contoura}]{
\includegraphics[scale=1]{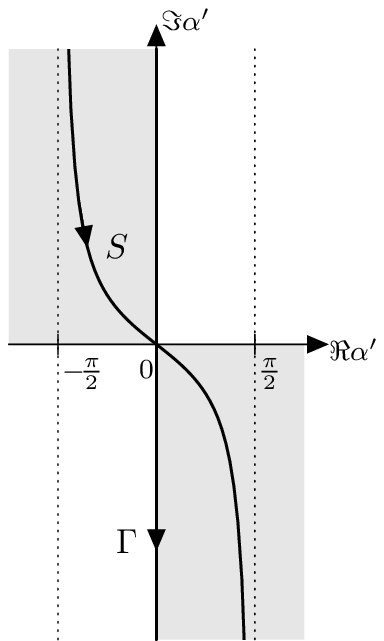}}
\hspace{-4mm}
\subfigure[\label{fig:contourb}]{
\includegraphics[scale=1]{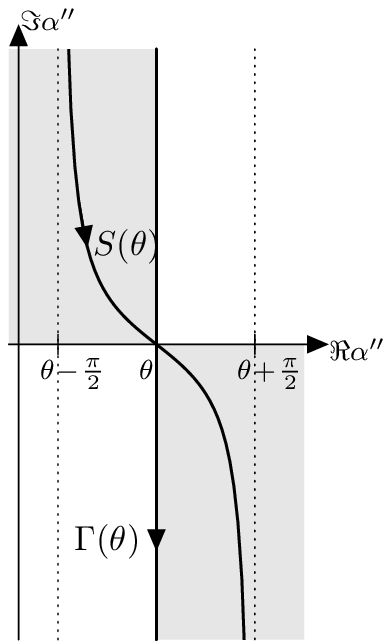}}
\caption{Integration contours for (a) Eqs.~\rf{eqn:pDdS} and \rf{eqn:pDdSImp}; (b) Eq.~\rf{eqn:pDdStheta}.}
\label{fig:Gamma}
\end{center}
\end{figure} 

One can obtain an expression suitable for far field (large $kr$) evaluation by deforming the contour in Eq.~\rf{eqn:DirSomm} onto the steepest descent contours $S_1$ and $S_2$ passing through the saddle points at $\alpha = \mp \pi$ (see Fig.~\ref{fig:gamma}). These contours are defined by $\Re{\alpha} = \mp \pi - \gd(\Im{\alpha})$ respectively, where $\gd(x)$ is the Gudermannian function with $\gd(x)=\cos^{-1}(\sech x)=\tan^{-1}(\sinh x)$, $0\leq \gd(x)\leq \frac{\pi}{2}$ for $x\geq0$, and $\gd(-x)=-\gd(x)$. To achieve the deformation one draws down the contour $\gamma_1$ onto $\gamma_2$, picking up residue contributions from any poles lying between $S_1$ and $S_2$. Assuming that $0\leq\theta\leq 3\pi/2$ and $0\leq\Re{\theta_0}\leq 3\pi/2$, one obtains\cite{Rawlins90}
\begin{align}
p^D(r,\theta) &= 
\He[\pi-|\theta-\Re{\theta_0}-\mathrm{gd}(\Im{\theta_0})|]\re^{-\ri k r \cos(\theta-\theta_0)}\notag\\
&\;\;-\He[\pi-|\theta+\Re{\theta_0}+\mathrm{gd}(\Im{\theta_0})|]\re^{-\ri k r \cos(\theta+\theta_0)}\notag\\
&\;\;-\He[\pi-|\theta+\Re{\theta_0}+\mathrm{gd}(\Im{\theta_0})-3\pi|]\re^{\ri k r \cos(\theta+\theta_0)}\notag\\
&\;\;+ p^D_d(r,\theta),
\label{eqn:SommTotal}
\end{align}
where
\begin{align}
\label{eqn:pDdS}
 p^D_d(r,\theta)=-\frac{1}{3\pi \ri } \int_{S}  \re^{\ri k r \cos{\alpha'}}\sin \frac{2}{3}\theta_0 \,G(\theta+\alpha',\theta_0)\,\rd \alpha',
\end{align}
\begin{align*}
\label{}
G(\alpha,\theta_0) &=  \frac{1}{\cos \frac{2}{3}\theta_0 - \cos\frac{2}{3}(\alpha-\pi)} \notag \\ &\;\;-  \frac{1}{\cos \frac{2}{3}\theta_0 - \cos\frac{2}{3}(\alpha+\pi)},%
\end{align*}
and the integrals over $S_1$ and $S_2$ have been combined into a single integral over the contour $S$ illustrated in Fig.~\ref{fig:contoura} by the changes of variable $\alpha = \alpha' \mp\pi$. 
Note that Eq.~\rf{eqn:SommTotal} corrects a sign in one of the exponents in the corresponding formula in Ref.~\cite{Rawlins90} (Eq.~(26) on p.~166). 
Note also that in Ref.~\cite{Rawlins90} Rawlins uses the equivalent representation 
\begin{align*}
\label{}
G(\alpha,\theta_0)=\frac{\sqrt{3}\sin \frac{2}{3}\alpha}{\left(\frac{1}{2}+\cos\frac{2}{3}(\alpha-\theta_0)\right)\left(\frac{1}{2}+\cos\frac{2}{3}(\alpha+\theta_0)\right)}.
\end{align*}

The first term in Eq.~\rf{eqn:SommTotal} represents the incident wave, and the second and third terms in Eq.~\rf{eqn:SommTotal} represent the two reflected waves. Note that the location of the zone boundary across which these waves are switched on/off by the Heaviside function prefactors is shifted compared to the case of purely real $\theta_0$. The shift agrees exactly with that derived in Ref.~\cite{Bertoni77} for the case of an inhomogeneous plane wave incident on a half plane. Note also that at a zone boundary the argument of one of the Heaviside functions in Eq.~\rf{eqn:SommTotal} equals zero, and a pole lies on the contour of integration. In this case a principal value integral should be taken in Eq.~\rf{eqn:pDdS} (for consistency with the assumption that $H[0]=1/2$); the same convention applies to the other decompositions in Eqs.~\rf{eqn:ptotal}, \rf{eqn:SommTotalEdgeSource}, \rf{eqn:SommTotalImp} and \rf{eqn:SommTotalImpEdgeSource}.

The edge source representation can then be obtained from Eq.~\rf{eqn:pDdS} as follows. First deform the contour of integration from $S$ to the imaginary axis $\Gamma$ in Fig.~\ref{fig:contoura} (equivalently, deform the original contours $\gamma_1$ and $\gamma_2$ onto the vertical contours $\Gamma_1$ and $\Gamma_2$ in Fig.~\ref{fig:gamma} rather than $S_1$ and $S_2$ before changing variables). Then write the resulting integral as an integral over the positive imaginary axis only, applying the identity%
\begin{align}
\label{eqn:GBetaIdentity}
\sin \frac{2}{3}\theta_0 \,\left(G(\theta+\ri \eta,\theta_0)+G(\theta-\ri \eta,\theta_0)\right) = -\beta^D,
\end{align}
which follows from the identity
\begin{align*}
&\frac{\sin a}{\cos a - \cos(b+c)} + \frac{\sin a}{\cos a - \cos(b-c)}\notag \\
&= - \frac{\sin (a+b)}{\cos c - \cos(a+b)}- \frac{\sin (a-b)}{\cos c - \cos(a-b)}.
\end{align*}
Finally, parametrise the resulting integral by 
\begin{align}
\label{eqn:ChangeVars}
 \alpha' = \ri\eta = \ri \sinh^{-1} (z/r),\qquad 0<z<\infty,
\end{align}
so that $r\cos\alpha'=r\cosh\eta = l$ and $d\alpha'/dz = \ri/l$, giving
\begin{align}
p^D(r,\theta) &= 
\He[\pi-|\theta-\Re{\theta_0}|]\re^{-\ri k r \cos(\theta-\theta_0)}\notag\\
&\quad-\He[\pi-|\theta+\Re{\theta_0}|]\re^{-\ri k r \cos(\theta+\theta_0)}\notag\\
&\quad-\He[\pi-|\theta+\Re{\theta_0}-3\pi|]\re^{\ri k r \cos(\theta+\theta_0)}\notag\\
&\quad+ p^D_{d,{\rm edge}}(r,\theta),
\label{eqn:SommTotalEdgeSource}
\end{align}
where $p^D_{d,{\rm edge}}(r,\theta)$ is given by the edge source integral in Eq.~\rf{eq:pwintegralspecial}. 

Clearly Eq.~\rf{eqn:SommTotalEdgeSource} reduces to Eq.~\rf{eqn:ptotal} when $0\leq \theta_0\leq  3\pi/2$. For complex $\theta_0$ there is a discrepancy in the arguments of the Heaviside functions between Eqs.~\rf{eqn:SommTotal} and \rf{eqn:SommTotalEdgeSource}; this is simply a consequence of the contour deformation from $S$ to $\Gamma$. Accordingly, Eq.~\rf{eqn:SommTotalEdgeSource} does not represent a decomposition of the field into `geometrical acoustics' and `diffracted' components, at least not in the sense understood in Ref.~\cite{Bertoni77}. Such discrepancies will have implications for the physical interpretation of the edge source integral derived for the impedance problem in the following section.

\section{Right-angled impedance wedge}
The main aim of this paper is to generalise Eq.~\rf{eqn:ptotal} to the case of impedance boundary conditions,
\begin{align}
\label{eqn:ImpBC}
\frac{\partial p}{\partial n} = \ri k \mu p,
\end{align}
where the unit normal vector $n$ points into the wedge, and $\mu$ is the complex admittance (inversely proportional to the impedance), which is assumed to satisfy $\Re{\mu}\geq0$, so as to prohibit energy creation at the boundary. 

It will be assumed that $\mu$ takes a constant value on each of the two wedge faces, but that the value of this constant may be different for the two wedge faces. Following Ref.~\cite{Rawlins90}, in order to simplify later formulas $\mu$ will be written
\begin{align}
\label{eqn:mudef}
\mu = 
\begin{cases}
-\sin{\theta_1}%
, & \textrm{on }\theta=0,\\
\;\;\;\cos{\theta_2}, & \textrm{on }\theta=3\pi/2,
\end{cases}
\end{align}
where $\theta_1,\theta_2$ are complex angles such that %
\begin{align*}
\label{}
\pi\leq \Re{\theta_1} \leq 3\pi/2,  \qquad %
0\leq \Re{\theta_2} \leq \pi/2. %
\end{align*}
Recalling that the wedge faces are given respectively by $y=0$, $x>0$ ($\theta=0$) and $x=0$, $y<0$ ($\theta=3\pi/2$), the boundary condition in Eq.~\rf{eqn:ImpBC} can then be stated as
\begin{align*}
\label{}
\begin{cases}
\dfrac{\partial p}{\partial y} = \ri k \sin{\theta_1}p, & \textrm{on }\theta=0,\\[3mm]
\dfrac{\partial p}{\partial x} = \ri k \cos{\theta_2}p, & \textrm{on }\theta=3\pi/2.
\end{cases}
\end{align*}
Note that when $\mu$ takes the same value on both faces the angles $\theta_1$ and $\theta_2$ are related by $\theta_2=3\pi/2-\theta_1$, since 
$-\sin{\theta_1}=\cos{(\frac{3\pi}{2}-\theta_1)}$. 

\subsection{Contour integral solution}
\label{sec:ImpContInt}
Rawlins shows that the total field $p^I$ for the impedance problem can be written as \cite{Rawlins90}
\begin{align}
\label{eqn:pIRep}
p^I(r,\theta) = L \,\tilde{p}^D(r,\theta),
\end{align}
where
\begin{align}
\label{eqn:LDef}
L &= \left(\frac{\partial^2}{\partial x \partial y}+ \ri k \left( \sin{\theta_1} \frac{\partial}{\partial x} + \cos{\theta_2} \frac{\partial}{\partial y}\right)\right. \notag \\ & \left. \phantom{\frac{1}{1}} \qquad\qquad - k^2 \sin{\theta_1}  \cos{\theta_2}\right)
\end{align}
and
\begin{align}
\tilde{p}^D(r,\theta)  &= \frac{\sin\frac{2}{3}\theta_{0}}{k^2\left(\cos\frac{2}{3}\theta_{2}-\cos\frac{2}{3}\theta_{1}\right)}\notag\\ &\quad \times \frac{1}{(\cos{\theta_2}-\cos{\theta_0})(\sin{\theta_0}-\sin{\theta_1})}\notag\\
& \quad\times
\sum_{j=0}^2 \frac{\cos\frac{2}{3}\theta_{j+2}-\cos\frac{2}{3}\theta_{j+1}}{\sin\frac{2}{3}\theta_{j}}p^{D,\theta_j},
\label{eqn:tildpDef}
\end{align}
where $\theta_{3}:=\theta_0$, $\theta_{4}:=\theta_1$ and $p^{D,\theta_j}$ denotes the Dirichlet solution for incident angle $\theta_j$ (this follows from Eqs.~(28)-(31) in Ref.~\cite{Rawlins90} combined with standard trigonometric identities). For completeness note that in Eq.~(29) on p.~167 of Ref.~\cite{Rawlins90}, $\cos \frac{2}{3}\theta_2$ in the numerator should be $\cos \frac{2}{3}\theta_1$. Rawlins' derivation of Eq.~\rf{eqn:pIRep} in Ref.~\cite{Rawlins90} is based on a trick first introduced by Williams  in Ref.~\cite{Williams65} to solve the analogous problem for a mixed Neumann/impedance wedge.

To evaluate the formula in Eq.~\rf{eqn:pIRep}, Rawlins uses the representation for $p^{D,\theta_j}$ given in Eq.~\rf{eqn:SommTotal}, but translates $\alpha'=\alpha''-\theta$ to give
\begin{align}
\label{eqn:pDdStheta}
 p^D_d(r,\theta)=-\frac{1}{3\pi \ri } \int_{S(\theta)}\!\!\!  \re^{\ri k r \cos{(\alpha''-\theta)}}\sin \frac{2}{3}\theta_0 \,G(\alpha'',\theta_0)\,\rd \alpha'',
\end{align}
where $S(\theta)$ is the translated version of $S$ passing through the saddle point at $\alpha''=\theta$ (see Fig.~\ref{fig:contourb}. This transformation greatly simplifies the application of the differential operator $L$ because the spatial variables $x$ and $y$ occur only in the exponential factor in the integrand in Eq.~\rf{eqn:pDdStheta}, with $r\cos(\alpha''-\theta) = x\cos\alpha'' + y \sin\alpha''$. Hence Eq.~\rf{eqn:pIRep} can be evaluated as\cite{Rawlins90}
\begin{align}
p^I(r,\theta) &= 
\He[\pi-|\theta-\theta_0|]\re^{-\ri k r \cos(\theta-\theta_0)}\notag\\
&\quad+R_1^I\He[\pi-\theta-\theta_0]\re^{-\ri k r \cos(\theta+\theta_0)}\notag\\
&\quad+R_2^I\He[\theta+\theta_0-2\pi]\re^{\ri k r \cos(\theta+\theta_0)}\notag\\
&\quad+T_1^I\He[\pi-\theta-\Re{\theta_1}-\gd(\Im{\theta_1})]\re^{-\ri k r \cos(\theta+\theta_1)}\notag\\
&\quad+T_2^I\He[\theta+\Re{\theta_2}+\gd(\Im{\theta_2})-2\pi]\re^{\ri k r \cos(\theta+\theta_2)}\notag\\
&\quad+ p^I_d(r,\theta),
\label{eqn:SommTotalImp}
\end{align}
where
\begin{align*}
\label{}
R_1^I &= \frac{\sin\theta_0+\sin\theta_1}{\sin\theta_0-\sin\theta_1},\\
R_2^I &= \frac{\cos\theta_0+\cos\theta_2}{\cos\theta_0-\cos\theta_2},\\
T_1^I &= \frac{2\sin \frac{2}{3}\theta_0(\cos \frac{2}{3}\theta_2 - \cos  \frac{2}{3}\theta_0)}{\sin \frac{2}{3}\theta_1(\cos \frac{2}{3}\theta_2 - \cos  \frac{2}{3}\theta_1)}
\\
& \quad \times
\frac{\sin\theta_1(\cos\theta_2-\cos\theta_1)}{(\sin\theta_0-\sin\theta_1)(\cos\theta_0-\cos\theta_2)}
,\\
T_2^I &= \frac{2\sin \frac{2}{3}\theta_0(\cos \frac{2}{3}\theta_1 - \cos  \frac{2}{3}\theta_0)}{\sin \frac{2}{3}\theta_2(\cos \frac{2}{3}\theta_2 - \cos  \frac{2}{3}\theta_1)}
\\
& \quad \times
\frac{\cos\theta_2( \sin\theta_2-\sin\theta_1)}{(\sin\theta_0-\sin\theta_1)(\cos\theta_0-\cos\theta_2)},
\end{align*}
and (changing variable back to $\alpha'$)
\begin{align}
\label{eqn:pDdSImp}
 p^I_d(r,\theta)&=-\frac{1}{3\pi \ri } \int_{S}  \re^{\ri k r \cos{\alpha'}}\frac{\sin \frac{2}{3}\theta_0}{\cos \frac{2}{3}\theta_2 - \cos  \frac{2}{3}\theta_1} \notag \\
 & \quad \times   \frac{(\sin(\alpha'+\theta)+\sin\theta_1)(\cos(\alpha'+\theta)+\cos\theta_2)}{(\sin\theta_0-\sin\theta_1)(\cos\theta_0-\cos\theta_2)} \notag \\
& \quad \times \sum_{j=0}^2 (\cos\tfrac{2}{3}\theta_{j+2}-\cos\tfrac{2}{3}\theta_{j+1})G(\alpha'+\theta,\theta_j)
\, \rd \alpha'.
\end{align}
The first term in Eq.~\rf{eqn:SommTotalImp} represents the incident field; the second and third terms represent reflected waves from the two wedge faces; the fourth and fifth terms represent surface waves propagating along the two wedge faces; and the final term $ p^I_d$ represents the diffracted field. 

A surface wave associated with the face $\theta=0$ is excited if $\pi\leq \Re\theta_1<\pi-\gd(\Im\theta_1)$, and in this case is confined to the angular region $0\leq\theta<\pi-\Re\theta_1-\gd(\Im\theta_1)$. Similarly, a surface wave associated with the face $\theta=3\pi/2$ is excited if $\pi/2-\gd(\Im\theta_2)<\Re\theta_2<\pi/2$, and in this case is confined to the angular region $2\pi-\Re\theta_2-\gd(\Im\theta_2)< \theta \leq 3\pi/2$. In terms of the admittance parameter $\mu$, recalling Eq.~\rf{eqn:mudef} one finds that surface waves are excited if
\begin{align*}
\label{}
\Re\mu\geq0, \qquad \Im\mu<\frac{\Re\mu}{\sqrt{1+(\Re\mu)^2}}.
\end{align*}
The surface waves (when they exist) decay exponentially with increasing distance both perpendicular to and along the face with which they are associated (unless $\mu$ is pure negative imaginary, in which case they maintain a constant amplitude along the face itself, decaying in the perpendicular direction).

The diffracted field $p^I_d$ can be approximated in the far field using the method of steepest descent, giving\cite{Rawlins90}
\begin{align}
\label{eqn:FarFieldImp}
 p^I_d(r,\theta) \sim D(\theta,\theta_0)\frac{\re^{\ri kr}}{\sqrt{kr}}+ O((kr)^{-3/2}), 
\end{align}
as $kr\to\infty$, where the diffraction coefficient 
\begin{align}
D(\theta,\theta_0) &= \frac{\sqrt{2}\re^{\ri\pi/4}}{3\sqrt{\pi}}\frac{\sin \frac{2}{3}\theta_0}{\cos \frac{2}{3}\theta_2 - \cos  \frac{2}{3}\theta_1} \notag \\
 & \quad \times   \frac{(\sin\theta+\sin\theta_1)(\cos\theta+\cos\theta_2)}{(\sin\theta_0-\sin\theta_1)(\cos\theta_0-\cos\theta_2)} \notag \\
& \quad \times \sum_{j=0}^2 (\cos\tfrac{2}{3}\theta_{j+2}-\cos\tfrac{2}{3}\theta_{j+1})\,G(\theta,\theta_j).
\label{eqn:DiffCoeffImp}
\end{align}
Note that Eq.~\rf{eqn:DiffCoeffImp} corrects a typographical error in Ref.~\cite{Rawlins90}: in Eq.~(36) of Ref.~\cite{Rawlins90}, $\sin{\theta}$ in the denominator should be $\sin{\theta_1}$. 
Note also that the approximation in Eq.~\rf{eqn:FarFieldImp} breaks down near zone boundaries (i.e., at values of $\theta$ for which the argument of one of the Heaviside functions in Eq.~\rf{eqn:SommTotalImp} equals zero). A more sophisticated far field approximation, valid uniformly across the zone boundaries, is given in Ref.~\cite{Rawlins09}.

\subsection{Edge source solution}
An edge source representation for $ p^I_d$ can be derived by closely following the procedure outlined in Section \ref{sec:DirDerivation} for the Dirichlet case. First deform the contour of integration from $S$ to $\Gamma$ (recall Fig.~\ref{fig:contoura}). Then write the resulting integral as an integral over the positive imaginary axis only. Simplifying the resulting expression requires slightly more work than in the Dirichlet case because the part of the integrand in Eq.~\rf{eqn:pDdSImp} not involving $G$ is no longer an even function of $\alpha'$ as it was in the Dirichlet case. To deal with this, first decompose
\begin{align*}
\label{}
(\sin(\alpha'+\theta)+\sin\theta_1)(\cos(\alpha'+\theta)+\cos\theta_2)\notag \\ \qquad\qquad = Q(\alpha') + \tilde Q(\alpha')
\end{align*}
into a sum of even and odd parts (with respect to $\alpha'$) 
\begin{align*}
\label{}
Q(\alpha') & =\cos\theta\sin\theta(\cos^2\alpha'-\sin^2\alpha') \notag \\
& \; + (\sin\theta_1\cos\theta + \cos\theta_2\sin\theta)\cos\alpha' + \sin\theta_1\cos\theta_2, \\
 \tilde Q(\alpha') & = (\cos^2\theta-\sin^2\theta)\cos\alpha'\sin\alpha'\notag \\
& \; +  (\cos\theta_2\cos\theta - \sin\theta_1\sin\theta)\sin\alpha'.
\end{align*}
Then deal with the contribution from the even part using Eq.~\rf{eqn:GBetaIdentity}, and that from the odd part using the identity
\begin{align}
\label{eqn:GBetaIdentityTilde}
\sin \frac{2}{3}\theta_0 \,\left(G(\theta+\ri \eta,\theta_0)-G(\theta-\ri \eta,\theta_0)\right) = -\ri\tilde\beta^D,
\end{align}
where
\begin{align*}
\label{}
\tilde\beta^{D} = \tilde\beta_1+ \tilde\beta_2-\tilde\beta_3-\tilde\beta_4,
\end{align*}
with 
\begin{align*}
\label{}
\tilde\beta_i = \frac{\sinh \frac{2}{3}\eta}{\cosh \frac{2}{3}\eta - \cos \frac{2}{3}\varphi_i}.
\end{align*}
Eq.~\rf{eqn:GBetaIdentityTilde} follows from the identity
\begin{align*}
&\frac{\sin a}{\cos a - \cos(b+c)} - \frac{\sin a}{\cos a - \cos(b-c)}\notag \\
&= \frac{\sin c}{\cos c - \cos(a+b)}- \frac{\sin c}{\cos c - \cos(a-b)}. 
\end{align*}
Finally, parametrising the contour using Eq.~\rf{eqn:ChangeVars} gives an edge source representation for the impedance solution:
\begin{align}
p^I(r,\theta) &= 
\He[\pi-|\theta-\theta_0|]\re^{-\ri k r \cos(\theta-\theta_0)}\notag\\
&\quad+R_1^I\He[\pi-\theta-\theta_0]\re^{-\ri k r \cos(\theta+\theta_0)}\notag\\
&\quad+R_2^I\He[\theta+\theta_0-2\pi]\re^{\ri k r \cos(\theta+\theta_0)}\notag\\
&\quad+T_1^I\He[\pi-\theta-\Re{\theta_1}]\re^{-\ri k r \cos(\theta+\theta_1)}\notag\\
&\quad+T_2^I\He[\theta+\Re{\theta_2}-2\pi]\re^{\ri k r \cos(\theta+\theta_2)}\notag\\
&\quad+ p^I_{d,{\rm edge}}(r,\theta),
\label{eqn:SommTotalImpEdgeSource}
\end{align}
where
\begin{align}
p_{d,{\rm edge}}^I &= -\frac{1}{3 \pi }   \int_{0}^{\infty}
\frac {\ {\rm e}^{\ri kl}} {l}\beta^I \, \rd z,
\label{eq:pwintegralspecialImp}
\end{align}
and
\begin{align*}
\label{}
\beta^I &= 
\frac{\sin\frac{2}{3}\theta_{0}}{\cos \frac{2}{3}\theta_2 - \cos  \frac{2}{3}\theta_1}\notag \\
& \quad \times
\frac{1}{(\sin\theta_0-\sin\theta_1)(\cos\theta_0-\cos\theta_2)}\notag \\
& \quad \times
\sum_{j=0}^2 \frac{\cos\frac{2}{3}\theta_{j+2}-\cos\frac{2}{3}\theta_{j+1}}{\sin\frac{2}{3}\theta_{j}}
\left(q\beta^{D,\theta_j}+\tilde q\tilde\beta^{D,\theta_j}\right), %
\end{align*}
where
\begin{align*}
\label{}
q & =\cos\theta\sin\theta(\cosh^2\eta+\sinh^2\eta) \notag \\
& \; + (\sin\theta_1\cos\theta + \cos\theta_2\sin\theta)\cosh\eta + \sin\theta_1\cos\theta_2, \\
 \tilde q & = (\sin^2\theta-\cos^2\theta)\cosh\eta\sinh\eta\notag \\
& \; +  (\sin\theta_1\sin\theta-\cos\theta_2\cos\theta)\sinh\eta,
\end{align*}
and, for $j=0,1,2$, $\beta^{D,\theta_j}$ and $\tilde\beta^{D,\theta_j}$ denote the functions $\beta^D$ and $\tilde\beta^D$ evaluated at incidence angle $\theta_j$. 

With regard to surface waves, note that the arguments of the Heaviside functions multiplying the fourth and fifth terms on the right-hand-side of Eq.~\rf{eqn:SommTotalImpEdgeSource} are always equal to zero except in the degenerate cases (i) $\Re\theta_1= \pi$ and $\theta=0$ and (ii) $\Re\theta_2=\frac{\pi}{2}$ and $\theta=\frac{3\pi}{2}$, respectively. Thus, recalling the discussion at the end of Section \ref{sec:ImpContInt}, one finds that if either $\theta_1$ or $\theta_2$ is such that surface waves are present, these surface waves must form part of the edge source integral in Eq.~\rf{eq:pwintegralspecialImp}. In this case, in the region where the surface waves exist it holds that
\begin{align}
\label{eqn:pdNonEqual}
p^I_{d,{\rm edge}}(r,\theta) \neq p^I_d(r,\theta) ,
\end{align}
so that the edge source integral in Eq.~\rf{eq:pwintegralspecialImp} cannot be associated solely with the diffracted field, as is the case for ideal (Dirichlet or Neumann) boundary conditions with a real incidence angle.  

Nonetheless, one can check that by applying the method of stationary phase to the integral in Eq.~\rf{eq:pwintegralspecialImp}, the far-field diffraction coefficient approximation in Eq.~\rf{eqn:FarFieldImp} is recovered. This does not contradict the above remarks (in particular Eq.~\rf{eqn:pdNonEqual}) since the surface wave contributions (when present) are exponentially small with respect to increasing $kr$, and hence are not picked up by the method of stationary phase.

It should be remarked that the existence of the edge source representation in Eq.~\rf{eqn:SommTotalImpEdgeSource} is of mainly theoretical interest (for example in the development of approximate solutions for finite edges and of edge integral equation formulations of scattering problems - see, e.g., Ref.~\cite{AsheimSvensson13}). For numerical computations of the infinite wedge solution at medium to high frequencies the expression in Eq.~\rf{eqn:SommTotalImp} should be used rather than that in Eq.~\rf{eqn:SommTotalImpEdgeSource}, because of the faster convergence of the integral in Eq.~\rf{eqn:pDdSImp} compared to that in Eq.~\rf{eq:pwintegralspecialImp}.
\section{Conclusions}\label{sect:conclusions}

A secondary edge source representation has been presented (in Eq.~\rf{eqn:SommTotalImpEdgeSource}) for the exact solution of scattering of a plane wave at perpendicular incidence on a right-angled impedance wedge. 
When the impedance parameters are such that surface waves are present, the edge source integral cannot be associated solely with the diffracted field, as in the case of ideal (Dirichlet or Neumann) boundary conditions, because it also incorporates the surface waves. 

A similar edge source representation should also be possible for general wedge angles, starting from the contour integral solutions in Refs.~\cite{Pierce78} and \cite{Osipov99}. But the analysis, and the resulting edge integral, are expected to be significantly more complicated than in the right-angled case considered here, for which one has the particularly simple contour integral solution provided by Ref.~\cite{Rawlins90}. 
Another interesting problem would be the derivation of edge source representations for more general incident waves, for example line source or point source excitation. However, to the present author's knowledge no convenient contour integral solution exists for these cases. Certainly the expressions obtained would be significantly more complicated than those obtained here for plane wave incidence; for a start, the Green's function for a line source or point source above an impedance boundary cannot be obtained by the method of images, as it can in the plane wave case. These generalisations are left for future work.
%

%
%
%
%
%
%
%
%
%
%
%

%
%
%
%
%

%
%
%

%
%
%
%

\end{document}